\documentclass[12pt,reqno]{amsart}

\usepackage[margin=0.6in]{geometry}
\usepackage{amsmath}
\usepackage{amsfonts}
\usepackage{amssymb}
\usepackage{amsthm}
\usepackage{mathtools}
\usepackage{graphicx,caption, color, hyperref}
\usepackage{parskip}
\usepackage{hyperref}
\usepackage{float}
\usepackage{multicol}

\numberwithin{equation}{section}
\newcounter{theorem}
\newtheorem{thm}[theorem]{Theorem}
\newtheorem{conj}[theorem]{Conjecture}
\newtheorem{cor}[theorem]{Corollary}
\newtheorem{prop}[theorem]{Proposition}
\newtheorem{lem}[theorem]{Lemma}

\theoremstyle{definition}
\newtheorem{defn}[theorem]{Definition}
\newtheorem{expl}[theorem]{Example}
\newtheorem{note}[theorem]{Note}

\numberwithin{theorem}{section}

\newcommand{\rog}{\text{Rog}(\mathbb{P})}

\newcommand{\set}[1]{\left\{#1\right\}}
\newcommand{\ord}{\text{ord}}

\title{A study of Cunningham bounds through Rogue Primes}

\author{A. Bhardwaj, L. Degen, R. Petkov, S. Stanbury.}

\begin{document}

\maketitle
\begin{abstract} If $p$ is prime, a sequence of prime numbers $\set{p, 2p+1, 4p+3,...,2^{n-1}(p+1)-1}$ is called a Cunningham chain. These are finite sequences of prime numbers, for which each element but the last is a Sophie Germain prime. It is conjectured that there are arbitrarily large such Cunningham chains, and these chains form an essential part of the study of Sophie Germain primes. In this paper,
we aim to significantly improve existing bounds for the length of Cunningham chains by considering their behaviour in the framework of what we will define as $\emph{rogueness}$. 
\end{abstract}

\section*{Introduction}
The congruence class of the base prime modulo $n$, where $n$ is prime, directly impacts the length of the subsequent Cunningham chain. This is particularly important when considering a potential strategy for finding arbitrarily long chains. Suppose $n = 5$, and we wish to generate a Cunningham chain of the first kind. Then if $a_{1}\equiv 1\mod 5$, we have
$$ a_2 \equiv 3 \mod 5\implies a_3 \equiv 2 \mod 5\implies a_4 \equiv 0 \mod 5,$$ forcing the chain to terminate at $a_3$. In fact, this demonstrates that any chain with length higher than 3 must have $a_1 \equiv 4 \mod 5$.  It follows similarly that any Cunningham chain of the second kind of length higher than 3 must have $a_1 \equiv 1 \mod 5$. The ability to place such a unique restriction allows us to bound chains for most initial elements. In addition, if we were able to find similar restrictions for other values of $n$, we could apply the Chinese Remainder Theorem on the base prime to significantly reduce the possibilities for chains of arbitrarily long length.
Moreover, for all $n$, a Cunningham chain of the first kind\footnote{see Definition \ref{defn:cunchain}} will never reach 0 modulo $n$ if $a_1 \equiv n-1 \mod n \equiv -1 \mod n$. This follows as $a_2 \equiv 2 \times -1 + 1 \equiv -1 \mod n$, with $a_k = -1 \mod n$ by induction. A similar result holds for a Cunningham chain of the second kind when $a_1 \equiv 1 \mod n$.

As seen for $n = 5$, we want as many congruence classes, such that a starting element from these classes would give a chain that would eventually hit $0$. Here, as seen above, for a chain of the $i^{\text{th}}$ kind, we ignore the equivalence class $(-1)^{i}$, and aim for every other equivalence class to fit our desired restriction. However, we observe that this is not possible for $n = 7$ by noting  
$$a_1 \equiv 2 \mod 7\implies a_2 \equiv 5 \mod 7\implies a_3 \equiv 4 \mod 7\implies a_4 \equiv 2 \equiv a_1 \mod 7,$$
leading to a repetition that will always be non-zero. This undesired loop for $n = 7$ does not allow us to enforce a unique restriction on the starting element, thus leading us to label $7$ as a $\emph{rogue}$ prime. 

In this paper, we will devote the first two sections to the construction of the idea of rogueness. The next sections apply these ideas to improve Cunningham bounds. The final two sections give further investigation into \emph{rogue loops}, and their relation to Cunningham chains. We also conjecture a logarithmic bound for these chains, for which we give numerical evidence. Finally, we conjecture a simple and deep result concerning divisors of the extensions of Cunningham chains, which is closely related to our logarithmic bound.

\section{Basic Definitions}
Throughout this section, we will be mainly concerned with odd prime numbers, and it is useful to have notation for the set of these primes.
\begin{defn}[Odd Prime Numbers]
We denote the set of $\textbf{odd prime numbers}$ to be
$$\mathbb{P}=\set{p\in\mathbb{N}\mid p>2,\gcd{(p,k)}=1,\forall k<p}.$$
\end{defn}
In this paper, we will focus on sequences of prime numbers. The basic pattern of the sequences we will study is  motivated by the notion of a Sophie Germain Prime.
\begin{defn}[Sophie Germain Prime]
A \textbf{Sophie Germain prime} is a prime number $p$, such that $2p+1$ is also prime.
\end{defn}
Remark that if for some odd prime $q<p$ we have $p\equiv\frac{q-1}{2}\mod q,$ then $p$ is not a Sophie Germain Prime. 

If $p$ is a Sophie Germain Prime, it is natural to wonder whether $2p+1$ is also a Sophie Germain prime, and if we can continue this process to get more primes. Whilst we keep hitting Sophie Germain Primes, we can continue this process, and in doing so we construct sequences of prime numbers, which are called Cunningham Chains.
\begin{defn}[Cunningham Chain]\label{defn:cunchain}
Let $p\in\mathbb{P}$. We define the $\textbf{Cunningham Chain}$ of the $i^{\text{th}}$ kind to be the finite sequence
$$a_{n}^{i}=2^{n-1}p+(-1)^{i+1}(2^{n-1}-1),\quad 1\leqslant n\leqslant k, a_{n}^{i}\in\mathbb{P}, a_{k+1}\not\in\mathbb{P}.$$
Here, $k$ depends on $i$ and $p$, $k=k(i,p)$ and the sequence is characterised by the recursion
$$a_{n+1}^{i}=2a_{n}^{i}+(-1)^{i+1}.$$
The prime $p$ is sometimes referred to as the \textbf{base} prime.
\end{defn}

For example, taking $p=2$, the Cunningham chain of the first/second kind would be (2, 5, 11, 23, 47) and (2, 3, 5) respectively.
It is believed that there exist infinitely many Cunningham chains, with the Bateman-Horn conjecture producing an asymptotic density estimate (\cite{loeh}). In fact, this would additionally allude to there being arbitrarily long chains, but as of August 2020, the longest known Cunningham chain is of length 19, discovered by Raanan Chermoni and Jaroslaw Wroblewski.
A $\textbf{bi-twin chain}$ is generated from a pair of twin primes $(n - 1, n + 1)$, as opposed to a single base prime. It takes the form $(n-1, n+1, 2n-1, 2n+1, ...)$ consisting of solely prime elements and can be considered as the intertwining of a Cunningham chain of the first kind with base prime $n - 1$ and a Cunningham chain of the second kind with base prime $n + 1$. 


It is useful to refer to the set of elements that form these Cunningham chains.
\begin{defn}[Cunningham Set]
Let $p\in\mathbb{P}$. Given the sequences $a_{n}^{i}$,
we define the $\textbf{Cunningham set of}$ 
\noindent $\textbf{ the }\mathbf{1^{\text{st}}/2^{\text{nd}}}\textbf{ kind}$ to be the set of elements in the Cunningham sequence. That is,
$$C_{i}[p]=\set{a_{k}^{i}\mid a_{k-1}^{i}, a_{k}^{i}\in\mathbb{P}:k\geqslant 2}\cup\set{a_{1}^{i}}, i=1,2.$$
\end{defn}
It is a simple consequence from group theory that if $p\in\mathbb{P}$, then the set $C_{i}[p]$ has size at most $p-1$. Moreover, from \cite{kanado} we can extract a better bound, namely for $p\geqslant 7$
\begin{equation}\label{eq:boundnotstrict}
    |C_{i}[p]|\leqslant\frac{p-3}{2},
\end{equation}
for $p\in\mathbb{P}$. In what follows, we
will consider a new approach for finding bounds for Cunningham Chains, by considering their concurrences, modulo primes.
\section{Rogue primes}
We now proceed to formalise the ideas mentioned in our introduction.
\begin{defn}[Rogue Set, Rogue Sequence]\label{defn:roguesetsequence}
Let $p$ be a prime. We define the $\textbf{rogue set of the}$ $\mathbf{1^{\text{st}}/2^{\text{nd}}}$ $\textbf{kind}$  by
$$A_{i}^{p}=\left(\mathbb{Z}/p\mathbb{Z}\right)^{\times}\setminus\set{(-1)^{i}\mod p},\quad i=1,2.$$
Let $g\in A_{i}^{p}$. Define the sequence
$$u_{n}^{i}=
\begin{cases}
g, n=1\\
2u_{n-1}^{i}+(-1)^{i+1}, n>1
\end{cases},$$
and set $g_{k}^{i}=\min\set{n\in\mathbb{N}\mid u_{n}^{i}=0}\cup\infty.$
We define the $\textbf{rogue sequence}$ of $g$ to be the sequence
$$\langle g\rangle_{i}^{p}=\set{u_{n}^{i}\mid n< g_{k}^{i}},$$
where here $k$ is as in Definition $\ref{defn:cunchain}$.

\end{defn}
To avoid having too much notation, when it will be clear, we will often denote the elements of a rogue sequence by $u_{n}^{i}$. We give a few examples of rogue sets and sequences to illustrate the above definition.
\begin{expl}\label{note:note1}
For $p\in\mathbb{P}$, the rogue sets of the $1^{\text{st}}$ and $2^{\text{nd}}$ kind are the sets (where here we make the abuse of notation by writing $x$ for $x \mod p$)
$$A_{1}^{p}=\set{1,2,3,...,p-2},\quad A_{2}^{p}=\set{2,3,...,p-2,p-1}.$$
Furthermore, for $p=5$ we have
$$\langle 1\rangle_{1}^{5}=\set{1,3,2}$$
whilst for $p=7$
$$\langle 2\rangle_{1}^{7}=\set{2,5,4,2,...},\quad\langle 2\rangle_{2}^{7}=\set{2,3,5,2,...}.$$
We thus have examples of both finite and infinite rogue sequences.
\end{expl}

We now give a formal definition for the primes that lack restriction on our Cunningham chains. To do so, we will first need to introduce terminology for these \emph{loops}, whose structure will generalise that of the rogue sequences $\langle 2\rangle_{1}^{7}$ and $\langle 2\rangle_{2}^{7}, $ seen in Example \ref{note:note1}. 
\begin{defn}[Rogue Loop, Rogue Prime of $1^{\text{st}} \text{ and  }2^{\text{nd}}$ kind]\label{defn:rogueprime}
Let $p\in\mathbb{P}$, and consider $A_{i}^{p}$. We call a periodic rogue sequence a $\textbf{rogue loop}$. We define $p$ to be $\textbf{rogue of the}$ $\mathbf{1^{\text{st}}/2^{\text{nd}}}$ $\textbf{kind}$ if $A_{i}^{p}$ contains a rogue loop. A prime which is not rogue is said to be a $\textbf{non-rogue prime}$.
\end{defn}
\begin{note}
    It is worth noting that in the case of rogue loops, we have that $g_{k}^{i}=\infty$.
\end{note}
It follows immediately by considering the closed form of the sequence $\langle(-1)^{i+1}\rangle_{i}^{p}$, that these always terminate, and hence are not rogue loops. Remark that for $p\in\mathbb{P}$, if we were to allow $(-1)^{i}$ to be in $A_{i}^{p}$, we would have as rogue sequences:
$$\langle (-1)^{i}\rangle_{i}^{p}=\set{(-1)^{i},(-1)^{i},(-1)^{i},...},$$
and hence we would always have a rogue loop, and all primes would be rogue.

Going back to Example \ref{note:note1}, we see that $\langle 1\rangle_{1}^{5}$ does not form rogue loop in $A_{1}^{5}$, whilst $\langle 2\rangle_{1}^{7}$ and $\langle 2\rangle_{2}^{7}$ do form rogue loops in $A_{1}^{7}$ and $A_{2}^{7}$ respectively. It follows that $7$ is a rogue prime of both the $1^{\text{st}}$ and $2^{\text{nd}}$ kind, whilst for $p=5$ we have
\begin{equation}\label{eq:eq1}
    \langle 1\rangle_{1}^{5}=\set{1,3,2},\langle2\rangle_{1}^{5}=\set{2},\langle3\rangle_{1}^{5}=\set{3,2},
\end{equation}
and hence $5$ is non-rogue of the first kind.
\begin{note}
Definition \ref{defn:rogueprime} is closely related to our motivational examples. If $p$ is non-rogue of the $1^{\text{st}}$ or $2^{\text{nd}}$ kind,
then the only criteria on a Cunningham generator $p_{1}$, of $1^{\text{st}} \text{ or  }2^{\text{nd}}$ kind is 
$$p_{1}\equiv-1 \mod p,\quad p_{1}\equiv 1\mod p.$$
\end{note}
In what follows, we will want to refer to the sets of rogue primes of the two kinds.
\begin{defn}[Set of Rogue Primes of $1^{\text{st}} \text{ and  }2^{\text{nd}}$ kind]
We define the $\textbf{set of rogue primes of the }\mathbf{1^{\text{st}}/2^{\text{nd}}}$ $\textbf{kind}$  to be
$$\text{Rog}(\mathbb{P})_{i}=\set{p\in\mathbb{P}\mid p\text{ is rogue of the } i^{\text{th}} \text{ kind}}.$$
\end{defn}

If we consider the rogue sets $A_{1}^{5}$ and $A_{2}^{5}$, we have
$$\langle 1\rangle_{1}^{5}=\set{1,3,2}\quad \langle -1\rangle_{2}^{5}=\set{4,2,3},$$
and hence every element of $A_{1}^{5}$ is contained in the sequence $\langle 1\rangle_{1}^{5}$, and every element of $A_{2}^{5}$ is contained in $\langle -1\rangle_{2}^{5}$. This leads to the idea of a rogue sequence being able to generate a rogue set.
\begin{defn}[Generated Rogue set, $p^{th}$ rogue generator]
Let $p\in\mathbb{N}$ be odd, consider the rogue set $A_{i}^{p}$. We say $A_{i}^{p}$ is a $\textbf{generated rogue set}$ if $\exists g\in A_{i}^{p}$ such that
$$A_{i}^{p}=\set{u_{n}^{i}\in\langle g\rangle_{i}^{p}}.$$
Note that here we ignore the ordering of the sequence, and it only matters that the elements of $A_{i}^{p}$ are the same as those in
$\langle g\rangle_{i}^{p}$.
We call such a $g$, a $\mathbf{p^{\text{th}}}\textbf{ rogue generator}$ of the $i^{\text{th}}$ kind.
\end{defn}
In $(\ref{eq:eq1})$, we found that the rogue sequence $\langle 1\rangle_{1}^{5}$ is an extension of the rogue sequence $\langle 3\rangle_{1}^{5}$, which itself was an extension of $\langle 2\rangle_{1}^{5}$. It is trivial that for any terminating rogue sequence of $\langle g\rangle_{i}^{p}$, the last element must be of the form
\begin{equation}\label{eq:queen}
    u_{g_{k}^{i}-1}\equiv \frac{p-(-1)^{i+1}}{2}\mod p.
\end{equation}
We can thus ask ourselves the question, how far back can we extend this sequence? In other words, which elements in $A_{i}^{p}$ can only ever be the first element of a rogue sequence. Clearly, for $p\in\mathbb{P}$, we have
$$\bigcup_{g\in A_{i}^{p}}\langle g\rangle_{i}^{p}=A_{i}^{p}.$$
We can further ask ourselves, what is the smallest subset $G\subset A_{i}^{p}$ such that
\begin{equation}\label{eq:liz2}
    \bigsqcup_{g\in G}\langle g\rangle_{i}^{p}=A_{i}^{p},
\end{equation}

 and if any such proper subsets of $A_{i}^{p}$ exist? Going back to Example \ref{note:note1}, in $A_{1}^{7}$ and $A_{2}^{7}$ we have
$$\langle 2\rangle_{1}^{7}=\set{2,5,4,2,...},\quad\langle 1\rangle_{1}^{7}=\set{1,3},\qquad \langle 2\rangle_{2}^{7}=\set{2,3,5,2,...},\quad\langle 6\rangle_{2}^{7}=\set{6,4},$$
and hence
$$A_{1}^{7}=\bigsqcup_{g\in\set{1,2}}\langle g\rangle_{1}^{7},\quad A_{2}^{7}=\bigsqcup _{g\in\set{2,6}}\langle g\rangle_{2}^{7}.$$

For now, we first consider this with $(\ref{eq:queen})$, and the following result becomes very immediate.
\begin{lem}\label{lem:lem40}
Let $p\in\mathbb{P}$. Then
$$A_{i}^{p} \text{ is generated }\iff\text{the } p^{\text{th}}\text{ rogue generator is given by } g_{i}\equiv (-1)^{i+1}\mod p.$$
\end{lem}
\begin{proof}
$(\Longrightarrow)$ Suppose $A_{i}^{p}$ is generated, then $\exists g\in A_{i}^{p}$ such that $\langle g\rangle_{i}^{p} = A_{i}^{p}.$ Fix $i=1$, and suppose $g\neq 1$. Then $1\in \langle g\rangle_{i}^{p} \setminus\set{g}$, and thus $\exists q\in \mathbb{N}$ such that
$$2q+1\equiv 1\mod p.$$
Thus
$$2q\equiv 0\mod p \implies q\equiv 0\mod p.$$
This is a contradiction since $0 \not\in A_{i}^{p}$. The proof is analogous for the case when $i=2$.
The converse follows by the definition of a rogue generator.  \end{proof}
Let $p\in\mathbb{P}$ and suppose for some $i\in\set{1,2}$, we have $p\not\in\rog_{i}$. It then follows that there are no rogue loops in $A_{i}^{p}$, and hence we can find a terminating sequence. In the same way as our motivation for Lemma \ref{lem:lem40}, we can extend such a sequence to get a generating sequence, and hence this raises the question as to which elements in $(\mathbb{Z}/p\mathbb{Z})^{\times}$ can be expressed as powers of $2$? This leads to the following remarkable result.
\begin{thm}\label{lem:lem20}
Let $p\in{\mathbb{P}}$. Then 
$$p\in\rog_{i}\iff 2\text{ is not a primitive root modulo } p.$$
\end{thm}
\begin{proof}
$(\Longrightarrow)$  Suppose $p\in\rog_{1}$, then there exists a rogue loop $\langle g\rangle_{1}^{p}$, given by
$$\langle g\rangle_{1}^{p}=\set{2^{n-1}(g+1)-1\mid n\in\mathbb{N}}.$$
By definition, $\forall n\in\mathbb{N}$ we have
$$2^{n-1}(g+1)\not\equiv 1\mod p,$$
and hence $2$ cannot be a primitive root modulo $p$. A similar process gives the same result for $p\in\rog_{2}.$

$(\Longleftarrow)$
Fix $i=1$ again, and suppose $p\not\in\rog_{1}$, then $\forall g\in A_{1}^{p}$, the rogue sequence $\langle g\rangle_{1}^{p}$ is finite. Hence for such $g, \exists K_{g}\in\mathbb{N}$ such that
$$2^{K_{g}-1}(g+1)-1\equiv0\mod p,$$
and hence
$$2^{K_{g}-1}\equiv (g+1)^{-1} \mod p.$$
Since $g$ was arbitrary and 
$$\langle 2\rangle\supseteq\set{(g+1)^{-1}\mid g\in A_{1}^{p}}=\set{g^{-1}\mid g\in (\mathbb{Z}/p\mathbb{Z})^{\times}\setminus\set{1}}=(\mathbb{Z}/p\mathbb{Z})^{\times}\setminus\set{1}.$$
It follows that $2$ is a primitive root modulo $p$. The proof is analogous for $i=2$.
\end{proof}
This leads to the following satisfying result, which will allow us to talk more freely about rogue primes.
\begin{cor}\label{lem:lem60}
Let $p\in\mathbb{P}.$ Then $p$ is a rogue prime of the $1^{\text{st}}$ kind $\iff$ $p$ is a rogue prime of the $2^{\text{nd}}$ kind.
\end{cor}
\begin{proof}
Follows from Theorem \ref{lem:lem20} since the right hand side is independent of $i$.
\end{proof}

\begin{defn}[Rogue Prime, Set of Rogue Primes]
Let $p\in\mathbb{P}$. Then $p$ is said to be a $\textbf{rogue prime}$ if $p$ is a rogue prime of the $1^{\text{st}} \text{ and/or  }2^{\text{nd}}$ kind.
We define the $\textbf{set of rogue primes}$ to be
$$\rog=\set{p\in\mathbb{P}\mid \text{$p$ is rogue}}.$$
\end{defn}
In other words, the rogue primes are simply the primes which do not have $2$ as a primitive root. Here is a list of the first $15$ of these:
$$\rog=\set{7,17,23,31,41,43,47,71,73,79,89,97,103,109,113,...}.$$
Suppose $A_{i}^{p}$ is generated, then any rogue sequence fits in as a subsequence of $\langle (-1)^{i+1}\rangle_{i}^{p}$ and hence cannot be a rogue loop. However, if $p\not\in\rog$, then every rogue sequence must terminate. The question remains, by how much can we extend these terminating sequences? We answer this with the following result.
\begin{lem}\label{lem:lem50}
Let $p\in\mathbb{P}$. Then
\begin{enumerate}
    \item $A_{i}^{p} \text{ is generated }\iff p\not\in\rog,$
    \item $A_{1}^{p} \text{ is generated }\iff A_{2}^{p} \text{ is generated.}$
\end{enumerate}
\end{lem}
\begin{proof}
\begin{enumerate}
\item $(\Longrightarrow)$ Let $p\in\rog$, and suppose $A_{i}^{p}$ is generated, and let $k\in A_{i}^{p}$. Then $\exists n\in\mathbb{N}$ such that
$$k=(-1)^{i+1}(2^{n}-1)\implies 2^{n}=(-1)^{i+1}k+1=(-1)^{i+1}(k+(-1)^{i+1}).$$
 Furthermore, since $p\in\rog$, by Lemma \ref{lem:lem20}, $2$ is not a primitive root modulo $p$, and hence $2$ has order at most $\frac{p-1}{2}$ in $(\mathbb{Z}/p\mathbb{Z})^{\times}$. Therefore
$$|\set{2^{n} \mod p\mid n\in\mathbb{N}}|\leqslant\frac{p-1}{2}.$$
Since $2$ is a primitive root modulo $3$, we can assume $p>3$. However
$$|\set{k+(-1)^{i+1}\mid k\in A_{i}^{p}}|=|(\mathbb{Z}/p\mathbb{Z})^{\times}\setminus\set{(-1)^{i+1}}|=p-2>\frac{p-1}{2},$$
giving a contradiction.

$(\Longleftarrow)$ Suppose $p\not\in\rog$ fix $i=1$, and consider
$$\langle 1\rangle_{1}^{p}=\set{2^{n}-1 \mod p\mid n\in\mathbb{N}}.$$ Then if $k\in A_{1}^{p},$ it follows that $k+1\in(\mathbb{Z}/p\mathbb{Z}_)^{\times}$, and hence $\exists N\in\mathbb{N}$ such that 
$$2^{N}=k+1\implies k=2^{N}-1,$$
and hence $A_{1}^{p}$ is generated. The analogous works for $i=2$.
\item Follows from $(1)$ immediately.
\end{enumerate}
\end{proof}
\begin{note}
We can thus partially answer the question raised in $(\ref{eq:liz2})$. When $p\not\in\rog$, we have that $A_{i}^{p}$ is generated, and by Lemma \ref{lem:lem40}, the generator is given by $(-1)^{i+1}$. Hence, in this case, it follows that $G=\set{(-1)^{i+1}}$, and we are done. The case for when $p\in\rog$, is more complex.
\end{note}
Having now worked through several results concerning rogue primes, rogue sets, and rogue sequences, we have built an intuition on the subject. We finish this section by noting that from Theorem \ref{lem:lem20}, it follows that Artin's Conjecture (\cite{murty}) for $n=2$ implies the existence of infinitely many rogue primes. This motivates further investigation of this notion of \emph{rogueness}.
\section{Cunningham Bounds}
\begin{figure}
\includegraphics[width=0.9\linewidth]{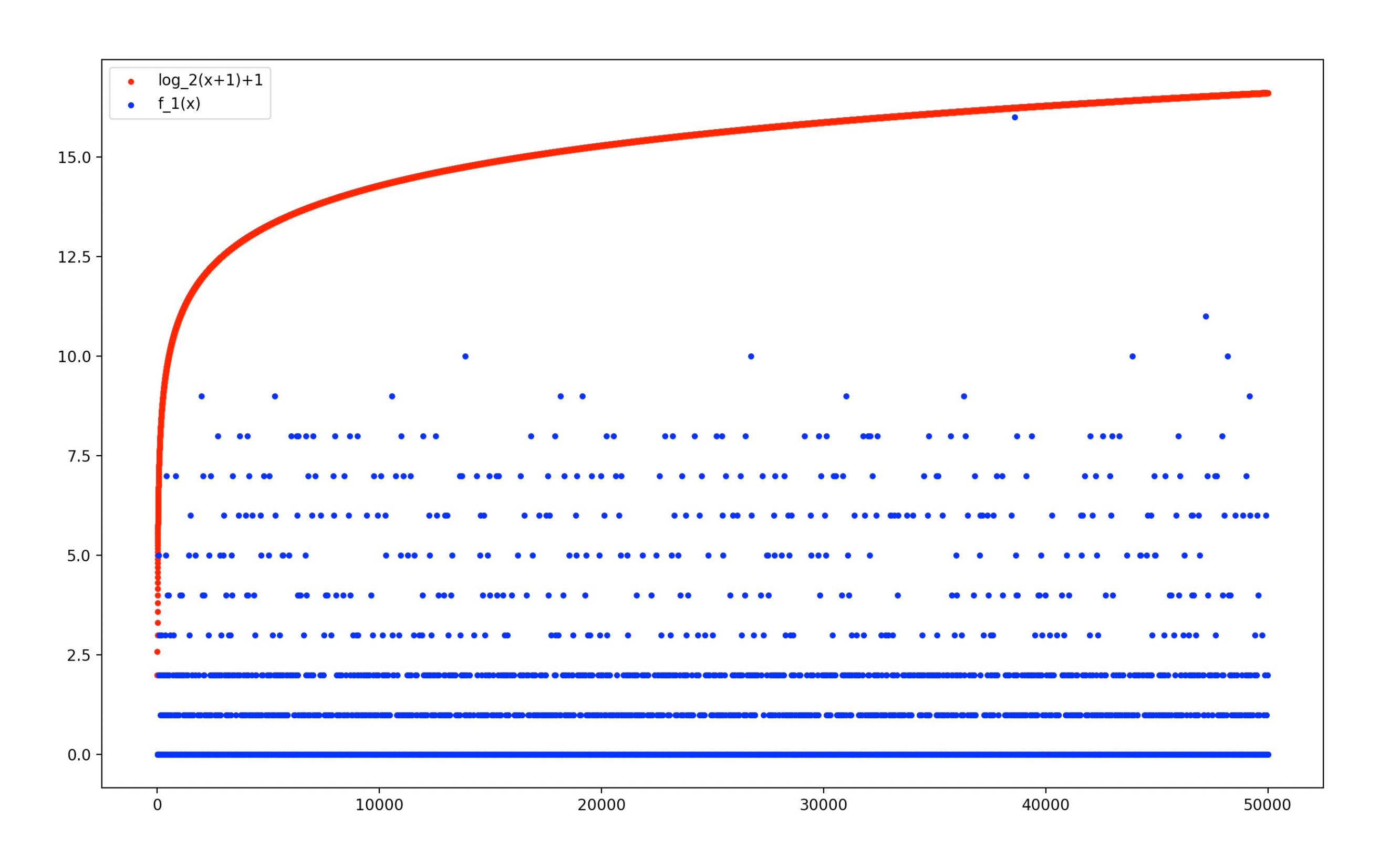}
\centering
\caption{$f_{1}(x)=\ord_{A_{1}^{k_{1}}}(x)-|C_{1}[x]|$}
\label{fig:bound1}
\end{figure}

\begin{figure}
\centering
\includegraphics[width=0.9\linewidth]{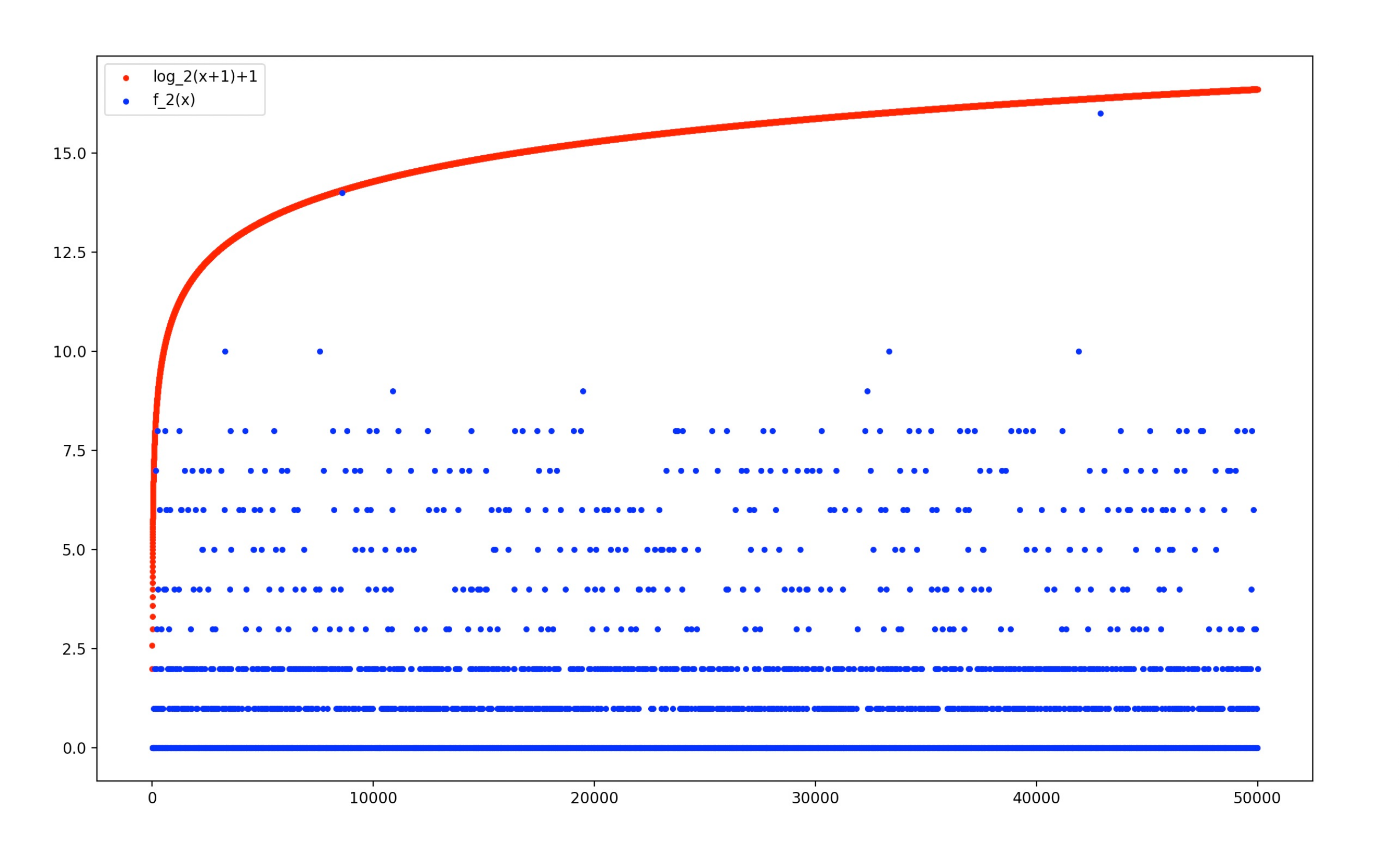}
\caption{$f_{2}(x)=\ord_{A_{2}^{k_{2}}}(x)-|C_{2}[x]|$}
\label{fig:bound2}
\end{figure}

In this section, we will apply the ideas developed in the previous section to study bounds for Cunningham chains; the main result of the section is Theorem \ref{thm:cunbound} which gives a new look at bounds for Cunningham chains, and $(\ref{eq:william})$ which gives a fast method of finding a bound for Cunningham chains. To do so, we will need a few more definitions, and to introduce terminology on the order of an element in a rogue sequence.
\begin{defn}\label{defn:q}
Let $n\in\mathbb{N}$, and let $i\in\set{1,2}$, we define the $Q_{i}^{n}$ to be the set
$$Q_{i}^{n}=\set{p\in\mathbb{P}\setminus\rog\mid p< n,  n\not\equiv (-1)^{i}\mod p}.$$
In other words, $Q_{i}^{n}$ is the set of non-rogue primes $p$ less than $n$, such that $n\in A_{i}^{p}$.
\end{defn}
We can express the elements of a generated rogue set by a rogue sequence. In the context of groups,  the order of an element vaguely corresponds to how far away that element is from the identity, in a multiplicative sense. We use this notion to give an analogous definition of order in such a generated rogue sense. This will essentially be the number of steps that element is from being divisible by the prime, upon which the rogue set is defined.
\begin{defn}[Order of an element in a generated rogue set]\label{defn:ordergenerated}
Let $p\in\mathbb{P}\setminus\rog$, and consider the generated set $A_{i}^{p}$. Let $q\in A_{i}^{p}$ for $i=1,2$. Then since $\langle (-1)^{i+1}\rangle_{i}^{p} = A_{i}^{p}$, $\exists k_{q}^{i}\leqslant p-2$ such that
$$q=u_{k_{q}^{i}}.$$
We define the $\textbf{order of q in}$ $A_{i}^{p}$ to be
\begin{equation}\label{eq:ord1}
    \ord_{A_{i}^{p}}(q)=p-1-k_{q}^{i}.
\end{equation}
\end{defn}
\begin{note}
It is worth noting that in the above definition, $k_{q}^{i}$ is the discrete logarithm base $2$ of $q+(-1)^{i}$. This is strongly related to the result in Theorem $\ref{lem:lem20}.$
\end{note}
We give some examples to illustrate how the order of elements in generated rogue sets is calculated.
\begin{expl}\label{expl:expl2}
In $A_{1}^{11}$ and $A_{2}^{11}$ we have
$$
\langle 1\rangle_{1}^{11}=\set{1,3,7,4,9,8,6,2,5},\quad
\langle 10\rangle_{2}^{11}=\set{10,8,4,7,2,3,5,9,6}.
$$
We thus have
$$\ord_{A_{1}^{11}}(7)=7,\quad,\ord_{A_{1}^{11}}(4)=6,\quad\ord_{A_{2}^{11}}(7)=6,\quad\ord_{A_{2}^{11}}(4)=7.$$
\end{expl}
\begin{note}
In Example \ref{expl:expl2}, we found that 
$$\ord_{A_{1}^{11}}(7)=\ord_{A_{2}^{11}}(4)=\ord_{A_{2}^{11}}(-7).$$
In fact there is nothing special about the fact that we are in $A_{i}^{11}$, and it follows for $p\not\in\rog$, $\forall k\in A_{i}^{p}\setminus\set{(-1)^{i+1}}$ we have from Definition \ref{defn:roguesetsequence}, and Lemma \ref{lem:lem40}
$$\ord_{A_{1}^{p}}(k)=\ord_{A_{2}^{p}}(-k).$$
\end{note}
We can now state and prove our main result of this section that will give us an improved bound for Cunningham Chains.
\begin{thm}[Rogue Bound for Cunningham Chains]\label{thm:cunbound}
Let $p\in\mathbb{P}$, and for $i\in\set{1,2},$ consider the Cunningham chain generated by $p$, then when $Q_{i}^{p}$ is non-empty we have
\begin{equation}\label{eq:godfamilyfootball}
    |C_{i}[p]|\leqslant\min\set{\ord_{A_{i}^{q}}(p)\mid q\in Q_{i}^{p}}.
\end{equation}
Moreover, in this case we have 
$$\min\set{\ord_{A_{i}^{q}}(p)\mid q\in Q_{i}^{p}}<p-1.$$
\end{thm}
\begin{note}\label{note:qipempty}
    The question of when $Q_{i}^{p}$ is non-empty is difficult. In Section \ref{sect:5}, we will develop terminology which will explain how Conjecture \ref{conj:conj1} allows us to bypass this difficulty. 
\end{note}
\begin{proof}
If $q\in Q_{i}^{p}$, then we have $q\not\in\rog$, and hence $A_{i}^{q}$ is generated. Furthermore, $q$ being in $Q_{i}^{p}$ implies that $p\in A_{i}^{q}$, and hence $\ord_{A_{i}^{q}}(p)$ is well defined. Thus $p\in \langle (-1)^{i+1}\rangle_{i}^{q}$, and $\exists k\leqslant q-2$ such that $p=u_{k}^{i}$, where here $k=q-1-\ord_{A_{i}^{q}}(p)$. Write $(a_{n})$ for the sequence of elements in $C_{i}[p]$, then  we have that the terms of the sequence $(a_{n})$ considered in $A_{i}^{q}$ are exactly $\langle u_{k}^{i}\rangle_{i}^{q}$. This is an equality of sequences. 
Since $u_{g_{k}^{i}}^{i}\equiv 0\mod q$, it follows
$$a_{\ord_{A_{i}^{q}}(p)}\equiv 0\mod q,$$
and thus $a_{\ord_{A_{i}^{q}}(p)}$ cannot be prime. Now, since $A_{i}^{q}$ is generated, it follows $\forall k\in A_{i}^{q}$
$$\ord_{A_{i}^{q}}(k)\leqslant q-1<p-1.$$
Hence  
$$|C_{i}[p]|\leqslant\min\set{\ord_{A_{i}^{q}}(p)\mid q\in Q_{i}^{p}}<p-1,$$
and the result follows.
\end{proof}
\begin{expl}\label{expl:expl1}
Take $p=89$, then 
$$C_{1}[89]=\set{89,179,359,719,1439,2879},$$
and hence $|C_{1}[89]|=6$. Using the previously known Cunningham bound, the best estimate we could have had was $$|C_{1}[89]|\leqslant 43.$$
We now consider our improved Cunningham bound.
$$Q_{1}^{89}=\set{11,13,19,29,37,53,59,61,67,83},$$
we find that
\begin{align*}
&\ord_{A_{1}^{11}}(89)=9,\quad\ord_{A_{1}^{13}}(89)=6,\quad\ord_{A_{1}^{19}}(89)=11,\quad\ord_{A_{1}^{29}}(89)=23,\quad\ord_{A_{1}^{37}}(89)=32,\\
&\ord_{A_{1}^{53}}(89)=22,\quad\ord_{A_{1}^{59}}(89)=9,\quad\ord_{A_{1}^{61}}(89)=25,\quad\ord_{A_{1}^{67}}(89)=38\quad\ord_{A_{1}^{83}}(89)=74,
\end{align*}
and hence 
$$|C_{1}[89]|\leqslant6.$$
In fact, in this case $|C_{1}[89]|=6,$ and hence our bound holds with equality. We will come back to this later on.
\end{expl}
\begin{note}
The bound given in Theorem \ref{thm:cunbound}, is quite a complicated one to conceptualise. However, since it is a minimum of a set, $\forall q\in Q_{i}^{p}$ it follows that
$$|C_{i}[p]|\leqslant\ord_{A_{i}^{q}}(p).$$
Hence, a more touchable upper bound would be
\begin{equation}\label{eq:william}
    |C_{i}[p]|\leqslant\ord_{A_{i}^{k_{i}}}(p),
\end{equation}
where $k_{i}=\min\set{q\in Q_{i}^{p}}.$ We thus only need to find the smallest non-rogue prime $k_{i}$ such that $p$ is not congruent to $(-1)^{i}$ modulo $k_{i}$, to find an upper bound for $C_{i}[p]$, that is better than $p-1$. For example, if $p\in\mathbb{P}$ such that $p>5$, and $p\not\equiv 4\mod 5$, then $5\in Q_{1}^{p}$, and hence $|C_{1}[p]|\leqslant 4.$ Similarly, if $p\not\equiv1\mod 5,$ then $5\in Q_{2}^{p}$, and hence $|C_{2}[p]|\leqslant4$. Moreover, we can even consider extreme cases, and find that this process is extremely fast: if we take $p=1122659$, then it follows that $|C_{1}[p]|=7$,
$$\min\set{q\in Q_{1}^{p}}=11,$$
and $\ord_{A_{1}^{11}}(p)=7$. We give a plot in Figures \ref{fig:bound1} and \ref{fig:bound2}, illustrating the accuracy of this more accessible bound. We refer the reader to \ref{app:f1} and \ref{app:f2} respectively, for the Python code of these figures.

Remark that the function $\log_{2}(x+1)+1$, is a bound for the error of each $f_{i}$, and hence we have
$$0\leqslant\ord_{A_{i}^{k_{i}}}(x)-|C_{i}[x]|<\log_{2}(x+1)+1,$$
giving the lower bound
$$|C_{i}[x]|>\ord_{A_{i}^{k_{i}}}(x)-\log_{2}(x+1)-1.$$
This bound however is only motivated from numerical calculations, and is more of an interesting observation for now.
\end{note}
\newpage

\section{Rogue Loops}
In $(\ref{eq:liz2})$, we discussed the partition of $A_{i}^{p}$ into disjoint rogue sequences. Whereas for $p\not\in\rog$, we found that the partition is a single rogue sequence. The case for $p\in\rog$ is more complicated, due to these rogue loops. We can thus ask ourselves; what length can these loops have, and how many distinct loops can there be? We begin to tackle these questions with a useful equivalence relation.
\begin{prop}\label{prop:rogueequiv}
Let $p\in\mathbb{P}$, the relation
$$R_{i}^{p}=\set{(a,b)\in A_{i}^{p}\times A_{i}^{p}\mid a\in\langle b\rangle_{i}^{p}\text{ or }b\in\langle a\rangle_{i}^{p}},$$
is an equivalence relation on $A_{i}^{p}$.
\end{prop}
\begin{proof}
Clearly $R_{i}^{p}$ is symmetric and reflexive, the only non-trivial condition is transitivity. Let $a,b,c\in A_{i}^{p}$, such that $(a,b),(b,c)\in R_{i}^{p}$.
Suppose $p\not\in\rog$, then 
$$\set{y\in\langle (-1)^{i+1}\rangle_{i}^{p}} =A_{i}^{p},$$
and hence $\exists k,l,m\in\mathbb{N}$ such that
$$a=u_{k}^{i},\quad b=u_{l}^{i},\quad c=u_{m}^{i}.$$
Without loss of generality, suppose $b\in \langle a\rangle_{i}^{p}$, and $c\in\langle b\rangle_{i}^{p}$. Then since 
$$\langle a\rangle_{i}^{p}=\set{u_{n}^{i}\mid n\geqslant k},\quad\langle b\rangle_{i}^{p}=\set{u_{n}^{i}\mid n\geqslant l},$$
it follows that $m\geq l\geq k$, and hence $c\in\langle a\rangle_{i}^{p}$, giving $(a,c)\in R_{i}^{p}$. Suppose that $p\in\rog$, and note that by Definition \ref{defn:rogueprime}, if $\langle g\rangle_{i}^{p}$ is a rogue loop in $A_{i}^{p}$, then $\forall k\in\langle g\rangle_{i}^{p}$, the sequence $\langle k\rangle_{i}^{p}$ is a shift of $\langle g\rangle_{i}^{p}$.
Thus if $(a,b)\in R_{i}^{p}$, and either one of $\langle a\rangle_{i}^{p}$, and $\langle b\rangle_{i}^{p}$ is a rogue loop, it follows that both rogue sequences are rogue loops. Hence if $(a,b)\in R_{i}^{p}$, then either both $\langle a\rangle_{i}^{p},\langle b\rangle_{i}^{p}$ are rogue loops, or both are terminating rogue sequences. Using the same logic for $(b,c)\in R_{i}^{p}$, we have that $\langle a\rangle_{i}^{p},\langle c\rangle_{i}^{p}$ are either both rogue loops or both terminating rogue sequences. If the former holds, then by our above reasoning it follows that $(a,c)\in R_{i}^{p}$. If the latter holds, then we can say without loss of generality that $a\in\langle b\rangle_{i}^{p}$, and $c\in\langle b\rangle_{i}^{p}$. It then follows $\exists m,k\in\mathbb{N}$ such that
$$a=v_{m}^{i},\quad c=v_{k}^{i},$$
where $\langle b\rangle_{i}^{p}=(v_{n}^{i})$. Suppose $k>m,$ it follows that
$$\langle a\rangle_{i}^{p}=\set{v_{n}^{i}\mid n\geqslant m}\ni c,$$
and hence $(a,c)\in R_{i}^{p}$. Other cases follow analogously.
\end{proof}
\begin{note}
We can thus partition $A_{i}^{p}$ into equivalence classes. It follows that the rogue loops are precisely the equivalence classes, except for the equivalence class of $(-1)^{i+1}$. Thus, $A_{i}^{p}$ can be written as the disjoint union of the rogue loops together with the elements of $\langle(-1)^{i+1}\rangle_{i}^{p}$. Of course, by Lemma \ref{lem:lem40}, we knew this for $p\not\in\rog$, but we have now covered all cases of $p\in\mathbb{P}$.
\end{note}
We now study these rogue loops to get a better understanding of this partition.
\begin{lem}\label{lem:terminatingsequences}
Let $p\in\rog$, and let $g\in A_{i}^{p}$, then
$$\langle g\rangle_{i}^{p}\text{ is a rogue loop}\iff g\not\in\langle(-1)^{i+1}\rangle_{i}^{p}.$$
\end{lem}
\begin{proof}
$(\Longrightarrow)$ Consider the rogue sequence
$$\langle(-1)^{i+1}\rangle_{i}^{p}=\set{u_{k}^{i}\mid k\in\mathbb{N}, k<g_{k}^{i}},$$
and suppose $g\in\langle(-1)^{i+}\rangle_{i}^{p}$. Then $\exists K\in\mathbb{N}$ such that $g=u_{K}^{i}$, and hence
$$\langle g\rangle_{i}^{p}=\set{u_{n}^{i}\mid n\geqslant K,n<g_{k}^{i}}.$$
It follows that $\langle g\rangle_{i}^{p}$ can not be periodic, and hence $\langle g\rangle_{i}^{p}$ is not a rogue loop.

$(\Longleftarrow)$ Suppose $\langle g\rangle_{i}^{p}$ is not a rogue loop. If $g=(-1)^{i+1}$, the result is obvious. Suppose $g\neq(-1)^{i+1}$. Then if 
$$\langle g\rangle_{i}^{p}=\set{v_{k}^{i}\mid k<g_{k}^{i}},$$
where $g_{k}^{i}=\min\set{n\in\mathbb{N}\mid v_{n}^{i}\equiv 0\mod p}$,
it follows by (\ref{eq:queen}) that 
$$v_{g_{k}^{i}-1}^{i}=\frac{p-(-1)^{i+1}}{2},$$
and hence $(g,\frac{p-(-1)^{i+1}}{2})\in R_{i}^{p}$, where $R_{i}^{p}$ is the equivalence relation on $A_{i}^{p}$ defined in Proposition \ref{prop:rogueequiv}. By a similar reasoning, $(\frac{p-(-1)^{i+1}}{2},(-1)^{i+1})\in R_{i}^{p}$, and since $R_{i}^{p}$ is an equivalence relation, it follows that $(g,(-1)^{i+1})\in R_{i}^{p}$. Hence $g\in\langle(-1)^{i+1}\rangle_{i}^{p}$ or $(-1)^{i+1}\in\langle g\rangle_{i}^{p}$. Suppose the latter holds, then $\exists N\in\mathbb{N}, N< g_{k}^{i}$ such that
$$2v_{N}^{i}+(-1)^{i+1}\equiv(-1)^{i+1}\mod p\implies v_{N}^{i}\equiv 0\mod p,$$
contradicting the minimality of $g_{k}^{i}.$ Hence $g\in\langle(-1)^{i+1}\rangle_{i}^{p}$, and the result follows.
\end{proof}
\begin{note}
    It follows from this result, that any terminating rogue sequence is a subsequence of the rogue sequence $\langle(-1)^{i+1}\rangle_{i}^{p}$. Of course, by Lemma \ref{lem:lem40}, we knew this for $p\not\in\rog$, and we can now use Lemma \ref{lem:terminatingsequences}, to generalise this for $p\in\rog$. 
\end{note}

We continue to study rogue loops. For a rogue prime $p$, and a rogue loop $\langle g\rangle_{i}^{p}$ in $A_{i}^{p}$, by the \textbf{length} of $\langle g\rangle_{i}^{p}$, we will mean the size of its underlying set.
We will denote the order of $2$ in $(\mathbb{Z}/p\mathbb{Z})^{\times}$ by $O_{p}(2)$. Take $p=7$, we have
$$A_{1}^{7}=\langle2\rangle_{1}^{7}\sqcup\langle 1\rangle_{1}^{7},\quad |\langle2\rangle_{1}^{7}|=3,\quad |\langle1\rangle_{1}^{7}|=2.$$
In this simple example, we notice that rogue loops have size one more than $\langle 1\rangle_{1}^{7}$. Furthermore, if $p=31$, we have
$$A_{1}^{31}=\langle2\rangle_{1}^{31}\sqcup\langle4\rangle_{1}^{31}\sqcup\langle27\rangle_{1}^{31}\sqcup\langle 29\rangle_{1}^{31}\sqcup\langle1\rangle_{1}^{31},$$
and a simple calculation gives that every rogue loop has size $5$, and that $|\langle1\rangle_{1}^{31}|=4$. Note that $2$ has order $3$ in $(\mathbb{Z}/7\mathbb{Z})^{\times}$, and has order $5$ in $(\mathbb{Z}/31\mathbb{Z})^{\times}.$ With this in mind, we provide the following powerful result.

\begin{thm}\label{thm:loops}
Let $p\in\rog$, then 
$$|\langle (-1)^{i+1}\rangle_{i}^{p}|=O_{p}(2)-1\leqslant \frac{p-3}{2}.$$
Furthermore, suppose that $\langle g\rangle_{i}^{p}$ is a rogue loop in $A_{i}^{p}$, then
\begin{equation}\label{eq:eq7}
    |\langle g\rangle_{i}^{p}|=O_{p}(2),
\end{equation}
and hence the maximum length of a rogue loop is $\frac{p-1}{2}$. Finally, the number of disjoint rogue loops in $A_{i}^{p}$ is given by
\begin{equation}\label{eq:eq6}
    \frac{p-1}{O_{p}(2)}-1.
\end{equation}
\end{thm}
\begin{note}
It is interesting to see that when $p\not\in\rog$, we have $O_{p}(2)=p-1$, and there are no rogue loops. Hence $(\ref{eq:eq6})$ holds for all $p\in\mathbb{P}$, and furthermore, for $p\not\in\rog$, it follows by Lemma \ref{lem:lem40}, that $$|\langle (-1)^{i+1}\rangle_{i}^{p}|=p-2=O_{p}(2)-1,$$
and hence the right equation in $(\ref{eq:eq7})$ also holds for all $p\in\mathbb{P}$.
\end{note}
\begin{proof}
Since $p\in\rog$, by Theorem \ref{lem:lem20}, $2$ is not a primitive root modulo $p$, and hence
\begin{equation}\label{eq:eq3}
    O_{p}(2)\leqslant\frac{p-1}{2}.
\end{equation}
Now let $\langle g\rangle_{i}^{p}$ be a rogue loop in $A_{i}^{p}$, then
$$\langle g\rangle_{i}^{p}=\set{2^{n-1}(g+(-1)^{i+1})-(-1)^{i+1} \mod p\mid n\in\mathbb{N}}.$$
Suppose $k,l\in\mathbb{N}$ such that 
$$2^{k-1}(g+(-1)^{i+1})-(-1)^{i+1}\equiv2^{l-1}(g+(-1)^{i+1})-(-1)^{i+1}\mod p,$$
it follows
$$2^{k}\equiv2^{l}\mod p\implies k\equiv l\mod O_{p}(2).$$
Hence the set $\set{2^{n-1}(g+(-1)^{i+1})-(-1)^{i+1} \mod p\mid n\in\mathbb{N}}$, contains $O_{p}(2)$ distinct values, and it follows
$$|\langle g\rangle_{i}^{p}|=O_{p}(2).$$
Furthermore, we have
$$\langle (-1)^{i+1}\rangle_{i}^{p}=\set{2^{n}(-1)^{i+1}-(-1)^{i+1}\mod p\mid n\in\mathbb{N}, n< g_{k}^{i}},$$
and since 
$$g_{k}^{i}=\min\set{n\in\mathbb{N}\mid 2^{n}(-1)^{i+1}-(-1)^{i+1}\equiv 0\mod p}=O_{p}(2),$$ 
it follows that $|\langle (-1)^{i+1}\rangle_{i}^{p}|=O_{p}(2)-1$. Then combining with $(\ref{eq:eq3})$, we get
$$|\langle(-1)^{i+1}\rangle_{i}^{p}|\leqslant\frac{p-1}{2}-1=\frac{p-3}{2}.$$
Finally by Proposition \ref{prop:rogueequiv}, we have
$$A_{i}^{p}=\bigsqcup_{[g]_{R}}\set{y\in\langle g\rangle_{i}^{p}}\implies p-2=|A_{i}^{p}|=\sum_{[g]_{R}}|\langle g\rangle_{i}^{p}|.$$
If we let 
$$X=\set{[g]_{R_{i}^{p}}\mid\langle g\rangle_{i}^{p}\text{ is a rogue loop}},$$
it follows that 
$$p-2=\sum_{[g]_{R_{i}^{p}}}|\langle g\rangle_{i}^{p}|=|\langle(-1)^{i+1}\rangle|+\sum_{[g]_{R_{i}^{p}}\in X}|\langle g\rangle_{i}^{p}|=O_{p}(2)-1+O_{p}(2)|X|,$$
and solving for $|X|$ gives the result.
\end{proof}
\begin{note}\label{note:note2}
Let $p\in\mathbb{P}$, and consider $C_{i}[p]$. It follows from the above result that if we can find a $q\in\rog$ such that $q<p$, and $p\in\langle(-1)^{i+1}\rangle_{i}^{q}$,
then
\begin{equation}\label{eq:eq4}
    |C_{i}[p]|\leqslant|\langle(-1)^{i+1}\rangle_{i}^{q}|\leqslant\frac{q-3}{2}<\frac{p-3}{2}.
\end{equation}
Hence, it would follow that the bound $(\ref{eq:boundnotstrict})$ is strict. In other words, if we can find a $q\in\rog$ such that $q<p$, and the congruence relation
\begin{equation}\label{eq:eq5}
    2^{n}\equiv(-1)^{i+1}p+1\mod q,
\end{equation}
has a solution, then $(\ref{eq:eq4})$ holds. 
\end{note}
We write this formally as a conjecture, which is numerically supported for the first 100000 rogue primes up to $N = 2157537$. We refer the reader to \ref{app:solns} for the Python code of this computation.
\begin{conj}\label{conj:solns}
Let $N\in\mathbb{N}$, let $i\in\set{1,2}$, then $\exists q\in\rog$ such that $N\in\langle(-1)^{i+1}\rangle_{i}^{q}$.
Hence $(\ref{eq:eq5})$ always has solutions.
\end{conj}
\section{Further investigation and Open Questions}\label{sect:5}
We begin this section by studying $(\ref{eq:godfamilyfootball})$. We have seen in Example \ref{expl:expl1}, that this equation can hold with equality and it is natural to ask ourselves whether this is always the case. If $p\in\mathbb{P}$, then one reason for $(\ref{eq:godfamilyfootball})$ to not hold with equality is if there exists a prime $q\in\rog$ such that $q<p,(p,(-1)^{i+1})\in R_{i}^{q}$, and
$$\ord_{A_{i}^{q}}(p)=|C_{i}[p]|.$$
Hence 
$2\cdot\max\set{C_{i}[p]}+(-1)^{i+1}$ could have $q$ as its only factor smaller than $p$, and thus $(\ref{eq:godfamilyfootball})$, would not hold with equality. This raises a further question, which is what if $2\max\set{C_{i}[p]}+(-1)^{i+1}$ has no divisors smaller than $p$? We tackle the first of these potential issues with an improved version of $Q_{i}^{p}$.
\begin{defn}[Corogueness]
Let $n\in\mathbb{N}$, a prime $q\in\mathbb{P}$ is \textbf{corogue} to $n$ if $(n,(-1)^{i+1})\in R_{i}^{q}.$ The set of primes corogue to $n$ is denoted
\begin{equation}\label{eq:corogue}
    \text{Rog}_{i}(n)=\set{p\in\mathbb{P}\mid p<n,\text{ and }(n,(-1)^{i+1})\in R_{i}^{p}}.
\end{equation}
\end{defn}
\begin{note}
    Clearly, $Q_{i}^{n}\subseteq\text{Rog}_{i}(n)$. In fact, $\text{Rog}_{i}(n)$ is simply the set $Q_{i}^{n}$ together with an appropriate set of rogue primes, depending on $n$.
\end{note}
We give an example to illustrate that this inclusion can be proper.
\begin{expl}

$$Q_{1}^{17}=\set{5,11,13},\quad\text{Rog}_{1}(17)=\set{5,7,11,13}.$$
\end{expl}
The set $\text{Rog}_{i}(n)$ is precisely the set we need in order to obtain an extended notion of order, which coincides with our original definition. The idea here is to use the results from Theorem \ref{thm:loops}: we can rewrite  $(\ref{eq:ord1})$ as 
$$\ord_{A_{i}^{p}}(q)=|\langle(-1)^{i+1}\rangle_{i}^{p}|-k_{q}^{i}=O_{p}(2)-1-k_{q}^{i}.$$
This motivates the following definition.
\begin{defn}[Order of an element in a corogue set]\label{defn:ordercorogue}
Let $n\in\mathbb{N},p\in\text{Rog}_{i}(n)$, 
we define the $\textbf{order of n}$ in $A_{i}^{p}$ to be
$$\ord_{A_{i}^{p}}(n)=O_{p}(2)-1-k_{p}^{i},$$
where 
$$u_{k_{p}^{i}}=n,$$
in the rogue sequence $\langle(-1)^{i+1}\rangle_{i}^{p}$.
\end{defn}
The only difference between Definition \ref{defn:ordercorogue} and Definition \ref{defn:ordergenerated}, is that for $n\in\mathbb{N}$, the number $\ord_{A_{i}^{p}}(n)$ is now defined for $p\in\text{Rog}_{i}(n)$, whereas previously, it was only defined for $p\in Q_{i}^{n}.$ We can thus give an improved version of Theorem \ref{thm:cunbound}, which is more likely to give equality. 
\begin{thm}\label{thm:refinement}
Let $p\in\mathbb{P}$, and suppose that $\text{Rog}_{i}(p)$ is non-empty. Then we have
\begin{equation}\label{eq:henry}
    |C_{i}[p]|\leqslant\min\set{\ord_{A_{i}^{q}}(p)\mid q\in\text{Rog}_{i}(p)}.
\end{equation}
\end{thm}
\begin{proof}
Repeat the argument of Theorem \ref{thm:cunbound}, using the extended terminology.
\end{proof}
To address the second issue we raised concerning equality of $(\ref{eq:godfamilyfootball})$, we  conjecture the existence of divisors of the first composite element of a rogue sequence, that are smaller then the base prime.
\begin{conj}\label{conj:conj1}
Let $p\in\mathbb{P}$, such that $p>3$, and let $k=2\max\set{C_{i}[p]}+(-1)^{i+1}$, then
$$\min\set{q\in\mathbb{P}: q\mid k}<p.$$
Hence $(\ref{eq:henry})$ holds with equality.
\end{conj}
The following note illustrates how Conjecture \ref{conj:conj1} justifies the choice of introducing rogue sets and rogue sequences to study Cunningham chains.
\begin{note}
    In Note \ref{note:qipempty} we discussed that the non-emptiness of $Q_{i}^{p}$ is difficult to study. The formulation of Theorem \ref{thm:refinement} translates this issue to the non-emptiness of $\text{Rog}_{i}(p)$. Conjecture \ref{conj:conj1} says that this is never the case. Hence, if Conjecture \ref{conj:conj1} holds true, in addition to saying that $(\ref{eq:henry})$ holds with equality, we can also remove the non-emptiness assumption on $\text{Rog}_{i}(p)$.
\end{note}
Another interesting consequence of Conjecture \ref{conj:conj1} is the following.
\begin{note}
If Conjecture \ref{conj:conj1} holds, then given $p\in\mathbb{P},\exists a\in\mathbb{P}$ such that $a<p$, and 
$$2^{|C_{i}[p]|}(p+(-1)^{i+1})+(-1)^{i}=a\cdot b,$$
for some $b\in\mathbb{N}$. It follows that
$$b=\frac{2^{|C_{i}[p]|}(p+(-1)^{i+1})+(-1)^{i}}{a}\geqslant\frac{2^{|C_{i}[p]|}(p+(-1)^{i+1})+(-1)^{i}}{p+(-1)^{i+1}}>2^{|C_{i}[p]|},$$
giving
$$|C_{i}[p]|<\log_{2}(b),$$
giving some kind of logarithmic bound.
\end{note}
With this remark, we give a simple logarithmic bound for a specific type of prime, namely the Mersenne Primes\footnote{See \cite{hill}, page 25,}.
\begin{lem}\label{lem:partialcunbound}
Let $p\in\mathbb{P}$, if $p=2^{n}-1$, for some $n\in\mathbb{N}$, then 
$$|C_{i}[p]|<\log_{2}(p+1).$$
\end{lem}
\begin{proof}
Suppose $p=2^{n}-1$, for some $n\in\mathbb{N}$, it follows that $\log_{2}(p+1)\in\mathbb{N}$, and
$$2^{\log_{2}(p+1)}(p+1)-1=(p+1)^{2}-1=p(p+2),$$
which is composite. Similarly, we have
$$2^{\log_{2}(p+1)}(p-1)+1=p^{2},$$
which is also composite.
\end{proof}
It is conjectured, that there are infinitely many such Mersenne primes, and hence we can imagine that this logarithmic bound could hold for infinitely many primes.
\begin{conj}\label{conj:logcunbound}
Let $p\in\mathbb{P}$, such that $p>5$, then we have
\begin{equation}\label{eq:logbound}
    |C_{i}[p]|<\log_{2}(p+1).
\end{equation}
\end{conj}
\begin{figure}[ht]
\includegraphics[width=0.9\linewidth]{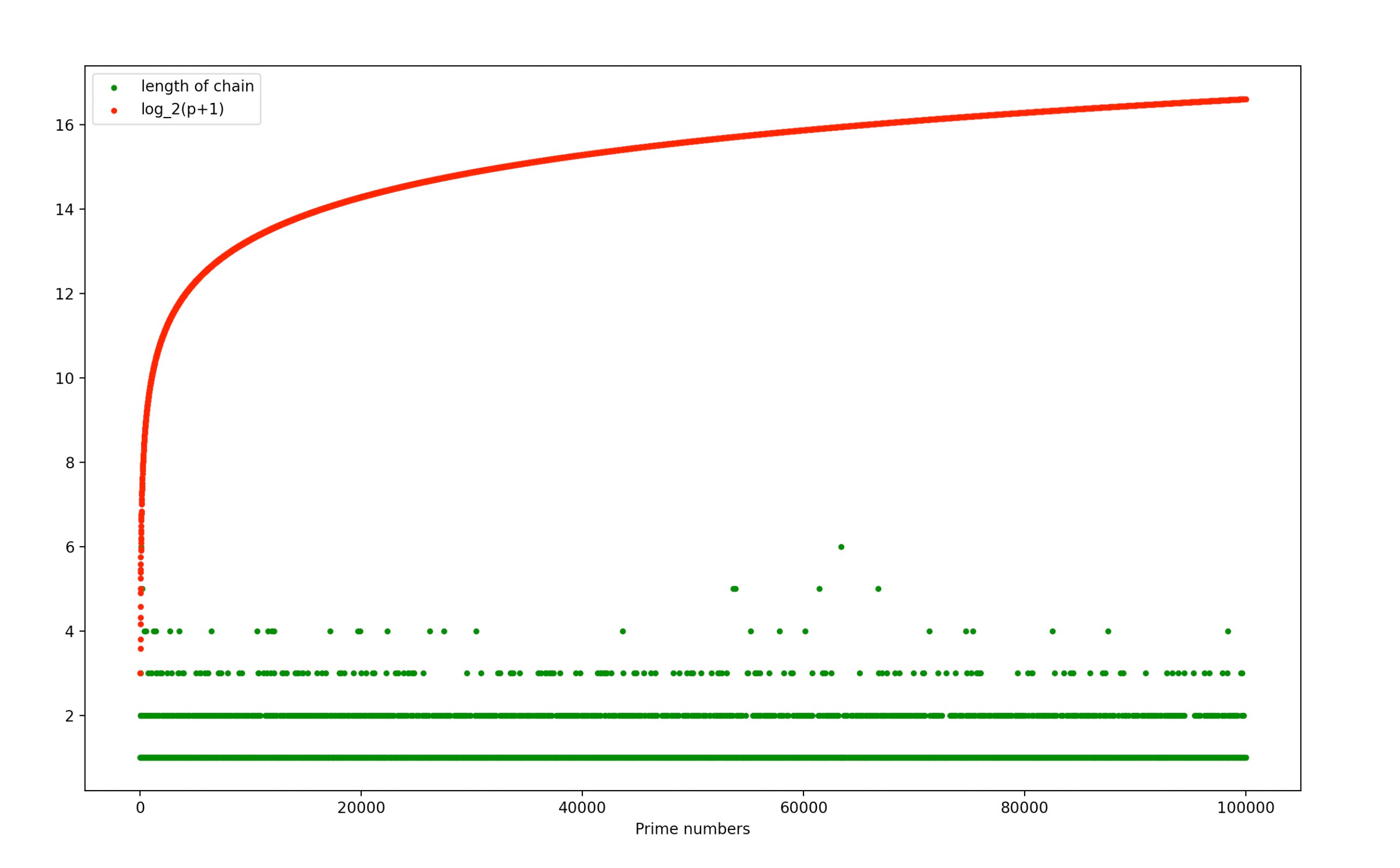}
\centering
\caption{Chains of the first kind}
\label{fig:bound3}
\end{figure}
\begin{figure}
    \centering
    \includegraphics[width=0.9\linewidth]{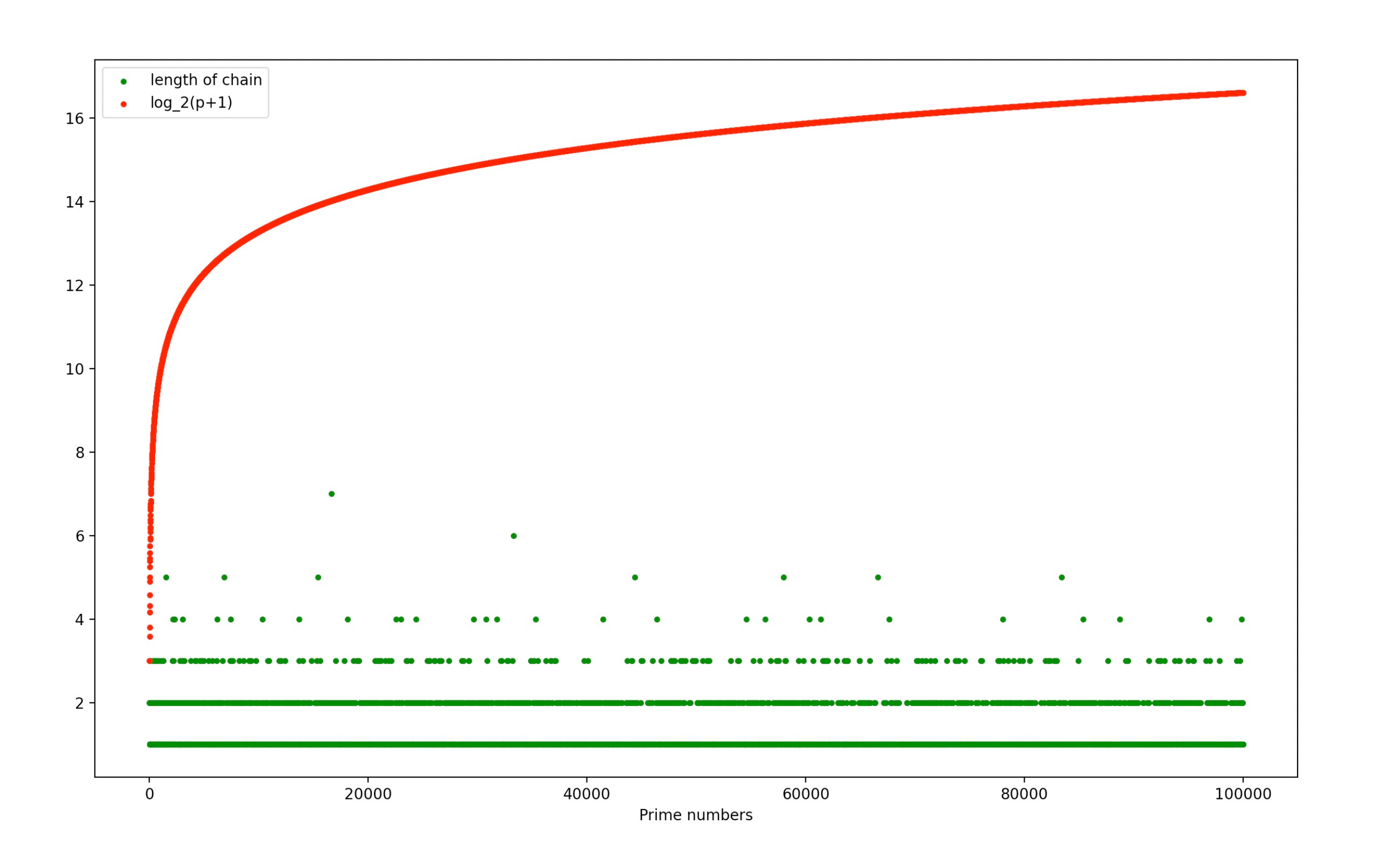}
    \caption{Chains of the second kind}
    \label{fig:bound4}
\end{figure}
We give graphs in Figures \ref{fig:bound3} and \ref{fig:bound4}, which illustrate this bound, and refer the reader to \ref{app:g1} and \ref{app:g2} respectively for the Python code of these figures.
\begin{thm}\label{thm:optimusprime}
Conjecture \ref{conj:logcunbound} $\implies$ Conjecture \ref{conj:conj1}.
\end{thm}
\begin{proof} We fix an $i\in\set{1,2}$, and show that by induction on the length of $|C_{i}[p]|$, by assuming  Conjecture \ref{conj:logcunbound}, we have Conjecture \ref{conj:conj1}. Suppose $p\in\mathbb{P}$ such that $|C_{i}[p]|=1$, then $p\neq 3$, and it follows that
$$p^{2}-(2p+1)>0,\quad p^{2}-(2p-1)>0.$$
Hence there exists a divisor of $2p+(-1)^{i+1}$ which is smaller than $p$.
Suppose $p\in\mathbb{P}$ such that $|C_{i}[p]|=n+1$, it follows that $|C_{i}[2p+1]|=n$, and hence there exists a prime divisor $q$ of 
$$k=2\max\set{C_{i}[p]}+(-1)^{i+1},$$
such that $q<2p+1$. Suppose for a contradiction, that $q\geqslant p$. Then we have
$$k\geqslant p^{2}\implies 2^{|C_{2}[p]|}(p-1)+1\geqslant p^{2}\implies|C_{2}[p]|\geqslant\log_{2}(p+1),$$
contradicting $(\ref{eq:logbound})$. 
\end{proof}
\begin{prop}
Let $p\in\mathbb{P}$, then $(\ref{eq:logbound})$ holds if and only if 
\begin{equation}\label{eq:cun2bound}
    2\max\set{C_{i}[p]}+(-1)^{i+1}\leqslant(p+1)^{2}.
\end{equation}
\end{prop}
\begin{proof}
$(\Longrightarrow)$ For $p>5$, the bound given in $(\ref{eq:logbound})$
gives $(\ref{eq:cun2bound})$. This follows by considering for the Cunningham Chain of the first kind
$$(p+1)^{2}-(2^{|C_{1}[p]|}(p+1)-1)>0\iff p>2^{|C_{1}[p]|-1}-1+\frac{1}{2}\sqrt{4^{|C_{1}[p]|}-4},
$$
and if $(\ref{eq:logbound})$ holds, then 
$$p>2^{|C_{1}[p]|}-1>2^{|C_{1}[p]|-1}-1+\frac{1}{2}\sqrt{4^{|C_{1}[p]|}-4}.$$
The analogous holds for chains of the second kind.

$(\Longleftarrow)$ In turn, if $(\ref{eq:cun2bound})$ holds it follows
$$p>2^{|C_{1}[p]|}-1>2^{|C_{1}[p]|-1}-1+\frac{1}{2}\sqrt{4^{|C_{1}[p]|}-4}>2^{|C_{1}[p]}-2,$$
giving
$|C_{1}[p]|\leqslant\log_{2}(p+2).$ 
Similarly, it follows that $|C_{2}[p]|\leqslant\log_{2}(p+2)$, and hence
$$|C_{i}[p]|\leqslant\log_{2}(p+2).$$
Then if $p+1=2^{n}$, for some $n\in\mathbb{N}$, it follows by Lemma \ref{lem:partialcunbound} that 
$|C_{i}[p]|<\log_{2}(p+1)$. If $p\neq 2^{n}$ for any $n\in\mathbb{N}$, then it follows that 
$$\lfloor\log_{2}(p+2)\rfloor=\lfloor\log_{2}(p+1)\rfloor<\log_{2}(p+1)<\log_{2}(p+2),$$
and hence 
$$|C_{i}[p]|<\log_{2}(p+1).$$
\end{proof}
We finish by remarking that Problem 5.1 in \cite{kanado}, implies that 
$$|C_{1}[p]|<\log_{2}(p)-1<\log_{2}(p+1),$$
giving Conjecture \ref{conj:logcunbound}. This last conjecture is heavily supported by numerical evidence, and hence gives further reason to investigate $(\ref{eq:logbound})$. By Theorem \ref{thm:optimusprime}, this would prove the very deep result that is Conjecture \ref{conj:conj1}, and give equality to $(\ref{eq:henry})$, giving an exact formula for the length of the Cunningham Chains.

\appendix

\section{Python Codes}

\subsection{Python code for $\mathbf{f_1(x):}$}\label{app:f1}
\begin{verbatim}
import matplotlib.pyplot as plt
import numpy as np
import math

def prime(n):
    for i in range(2, int(n / 2) + 1):
        if (n % i == 0):
            return 0
    return 1

def getSophieGermainPrime(startLimit, endLimit):
    l1 = []
    p, r = startLimit, endLimit
    if p == 1:
        p = 2
    for a in range(p, r + 1):
        k = 0
        for i in range(2, int(a / 2) + 1):
            if (a % i == 0):
                k = k + 1
                break
        if (k <= 0):
            x = prime(2 * a + 1)
            if (x == 1):
                l1.append(a)
    return (l1)

def print_C_len(p0):
    i = 0;
    list = []
    while (True):
        flag = 1;
        x = pow(2, i);
        p1 = x * p0 + (x - 1);

        for k in range(2, p1):
            if (p1 % k == 0):
                flag = 0;
                break;

        if (flag == 0):
            break;
        list.append(p1)
        i += 1;
    return (len(list))

def print_C_list(p0):
    i = 0;
    list = []
    while (True):
        flag = 1;
        x = pow(2, i);
        p1 = x * p0 + (x - 1);

        for k in range(2, p1):
            if (p1 % k == 0):
                flag = 0;
                break;

        if (flag == 0):
            break;
        list.append(p1)
        i += 1;
    return (list)

List_Primes=[]

for i in range(6,50000):
    if prime(i)==True:
        List_Primes.append(i)

Non_rogue_list = [set of Non rogue prime numers]
p = 0
SG_primes = getSophieGermainPrime(6, 50000)
Chain_length = []

x1 = [i for i in List_Primes]
for j in x1:
    p = print_C_len(j)
    Chain_length.append(p)

Q = []
List_of_k = []
for i in List_Primes:
    for j in Non_rogue_list:
        if j < i:
            Q.append(j)
        for m in Q:
            if (i + 1) % m == 0:
                Q.remove(m)
        if Q:
            k_i = min(Q)
            List_of_k.append(k_i)
            Q = []
            break
Order_list=[]
Chain_list=[]

index=0
for i in range(len(List_Primes)):
    new_k=List_of_k[i]
    remainder=List_Primes[i]%new_k
    while remainder !=0:
        remainder=remainder*2+1
        remainder=remainder%new_k
        index=index+1
    Order_list.append(index)
    index=0

array1 = np.array(Order_list)
array2 = np.array(Chain_length)
subtracted_array = np.subtract(array1, array2)
subtracted = list(subtracted_array)
index2=0

x=range(1,50000)
y=[]
for i in x:
    y.append(math.log(i+1,2)+1)
plt.scatter(x, y, marker='.', label="log_2{x+1}+1",c='r')
plt.scatter(List_Primes, subtracted, marker='.', label="f_1(x)",c='b')

plt.legend()
plt.show()

\end{verbatim}
\subsection{Python code for $\mathbf{f_2(x):}$}\label{app:f2}
\begin{verbatim}
import matplotlib.pyplot as plt
import numpy as np
import math

def prime(n):
    for i in range(2, int(n / 2) + 1):
        if (n % i == 0):
            return 0
    return 1

def getSophieGermainPrime(startLimit, endLimit):
    l1 = []
    p, r = startLimit, endLimit
    if p == 1:
        p = 2
    for a in range(p, r + 1):
        k = 0
        for i in range(2, int(a / 2) + 1):
            if (a % i == 0):
                k = k + 1
                break
        if (k <= 0):
            x = prime(2 * a - 1)
            if (x == 1):
                l1.append(a)
    return (l1)

def print_C_len(p0):
    i = 0;
    list = []
    while (True):
        flag = 1;
        x = pow(2, i);
        p1 = x * p0 - (x - 1);

        for k in range(2, p1):
            if (p1 % k == 0):
                flag = 0;
                break;

        if (flag == 0):
            break;
        list.append(p1)
        i += 1;
    return (len(list))

def print_C_list(p0):
    i = 0;
    list = []
    while (True):
        flag = 1;
        x = pow(2, i);
        p1 = x * p0 - (x - 1);

        for k in range(2, p1):
            if (p1 % k == 0):
                flag = 0;
                break;

        if (flag == 0):
            break;
        list.append(p1)
        i += 1;
    return (list)

List_Primes=[]

for i in range(6,50000):
    if prime(i)==True:
        List_Primes.append(i)

Non_rogue_list = [set of Non rogue prime numers]
p = 0
SG_primes = getSophieGermainPrime(6, 50000)
Chain_length = []

x1 = [i for i in List_Primes]
for j in x1:
    p = print_C_len(j)
    Chain_length.append(p)

Q = []
List_of_k = []
for i in List_Primes:
    for j in Non_rogue_list:
        if j < i:
            Q.append(j)
        for m in Q:
            if (i - 1) % m == 0:
                Q.remove(m)
        if Q:
            k_i = min(Q)
            List_of_k.append(k_i)
            Q = []
            break
Order_list=[]
Chain_list=[]

index=0
for i in range(len(List_Primes)):
    new_k=List_of_k[i]
    remainder=List_Primes[i]%new_k
    while remainder !=0:
        remainder=remainder*2-1
        remainder=remainder%new_k
        index=index+1
    Order_list.append(index)
    index=0

array1 = np.array(Order_list)
array2 = np.array(Chain_length)
subtracted_array = np.subtract(array1, array2)
subtracted = list(subtracted_array)
index2=0

x=range(1,50000)
y=[]
for i in x:
    y.append(math.log(i+1,2)+1)
plt.scatter(x, y, marker='.', label="log_2{x+1}+1",c='r')
plt.scatter(List_Primes, subtracted, marker='.', label="f_2(x)",c='b')

plt.legend()
plt.show()


\end{verbatim}
\subsection{Python code for $\mathbf{g_1(x):}$}\label{app:g1}
\begin{verbatim}
import matplotlib.pyplot as plt
import math
def prime(n):
    for i in range(2,int(n/2)+1):
        if(n%i==0):
            return 0
    return 1
def getSophieGermainPrime(startLimit,endLimit):
    l1=[]
    p,r=startLimit,endLimit
    if p==1:
        p=2
    for a in range(p,r+1):
        k=0
        for i in range(2,int(a/2)+1):
            if(a%i==0):
                k=k+1
                break
        if(k<=0):
            x=prime(2*a+1)
            if(x==1):
                l1.append(a)
    return(l1)

def print_C(p0):
    i = 0;
    list=[]
    while (True):
        flag = 1;
        x = pow(2, i);
        p1 = x * p0 + (x - 1);

        for k in range(2, p1):
            if (p1 % k == 0):
                flag = 0;
                break;

        if (flag == 0):
            break;
        list.append(p1)
        i += 1;
    return(len(list))
p = 0
SG_primes = getSophieGermainPrime(1,1000)

List_Primes=[]
for i in range(6,100000):
    if prime(i)==True:
        List_Primes.append(i)
y1=[]
x1 = [i for i in List_Primes]
for j in x1:
    p = print_C(j)
    y1.append(p)

y=[]
for i in List_Primes:
    y.append(math.log(i+1,2))

plt.scatter(x1, y1, marker='.', label="length of chain",c='g')
plt.scatter(List_Primes, y, marker='.', label="log_2(p+1)",c='r')
plt.xlabel('Prime numbers')
plt.legend()
plt.show()

\end{verbatim}
\subsection{Python code for $\mathbf{g_2(x):}$}\label{app:g2}
\begin{verbatim}
import matplotlib.pyplot as plt
import math
def prime(n):
    for i in range(2,int(n/2)+1):
        if(n%i==0):
            return 0
    return 1
def getSophieGermainPrime(startLimit,endLimit):
    l1=[]
    p,r=startLimit,endLimit
    if p==1:
        p=2
    for a in range(p,r+1):
        k=0
        for i in range(2,int(a/2)+1):
            if(a%i==0):
                k=k+1
                break
        if(k<=0):
            x=prime(2*a-1)
            if(x==1):
                l1.append(a)
    return(l1)

def print_C(p0):
    i = 0;
    list=[]
    while (True):
        flag = 1;
        x = pow(2, i);
        p1 = x * p0 - (x - 1);

        for k in range(2, p1):
            if (p1 % k == 0):
                flag = 0;
                break;

        if (flag == 0):
            break;
        list.append(p1)
        i += 1;
    return(len(list))
p = 0
SG_primes = getSophieGermainPrime(1,1000)

List_Primes=[]
for i in range(6,100000):
    if prime(i)==True:
        List_Primes.append(i)
y1=[]
x1 = [i for i in List_Primes]
for j in x1:
    p = print_C(j)
    y1.append(p)

y=[]
for i in List_Primes:
    y.append(math.log(i+1,2))

plt.scatter(x1, y1, marker='.', label="length of chain",c='g')
plt.scatter(List_Primes, y, marker='.', label="log_2(p+1)",c='r')
plt.xlabel('Prime numbers')
plt.legend()
plt.show()

\end{verbatim}

\subsection{Python code for Conjecture \ref{conj:solns} }\label{app:solns}
\begin{verbatim}
from itertools import islice
from sympy import nextprime, is_primitive_root

def generator():
    p = 2
    while(p:=nextprime(p)):
        if not(is_primitive_root(2, p)):
            yield p

roguePrimes = list(islice(generator(), 100000)) #OEIS A001122 Chai Wah Wu

c = max(roguePrimes)
#print(c)
valid = set() #numbers which are generated by plus (minus) 1


for prime in roguePrimes:
    generated = []
    a = 1
    while True:
        generated.append(a)
        a = (2*a + 1)%prime
        if a == 0:
            break
    valid.update(generated)

m = list(valid)

for i in range(1, c+1): #largest number+1 for which we can guarantee q
    if m[i-1] != i:
        print(i)
        break

\end{verbatim}

\section*{Acknowledgments} We are extremely thankful to Dr Keenan Kidwell for his incredible support and insightful feedback that has culminated in the paper as it stands. We would like to thank Melissa Ozturk for her contributions and insights during our group project in the final term of the 2021/22 academic year, which led to the writing of this paper. We are deeply thankful to Elliot Cocks for his calculations which helped us get these ideas started. 

\bibliographystyle{plain}
\bibliography{year2main.bib}

{\footnotesize  
\vspace*{0.5mm} 
 
 \begin{multicols}{2}

\noindent {\it Anand Bhardwaj}\\  
University College London\\
25 Gordon St, WC1H 0AY\\
London, United Kingdom\\
E-mail: {\tt anand.bhardwaj.20@ucl.ac.uk}\\ \\ 
\noindent {\it Luisa Degen}\\  
University College London\\
25 Gordon St, WC1H 0AY\\
London, United Kingdom\\
E-mail: {\tt luisa.degen.20@ucl.ac.uk}\\ \\ 

\noindent {\it Radostin Petkov}\\  
University College London\\
25 Gordon St, WC1H 0AY\\
London, United Kingdom\\
E-mail: {\tt radostin.petkov.20@ucl.ac.uk}\\ \\ 
\noindent {\it Sidney Stanbury}\\  
University College London\\
25 Gordon St, WC1H 0AY\\
London, United Kingdom\\
E-mail: {\tt sidney.stanbury.20@ucl.ac.uk}\\ \\

\end{multicols}

\end{document}